\documentclass[fleqn,leqno]{article}
\usepackage{amsfonts}

\newtheorem{theorem}{Theorem}

\input{tcilatex}
\begin{document}

\begin{center}
\textbf{\Large Asymptotic series and inequalities associated to some
expressions involving the volume of the unit ball}\bigskip 

\textbf{Cristinel Mortici$^{1,2}$}

\bigskip

$^{1}$Department of Mathematics, Valahia University of T\^{a}rgovi\c{s}te,%
\\[0pt]
Bd. Unirii 18, 130082 T\^{a}rgovi\c{s}te, Romania\\[0pt]
$^{2}$Academy of Romanian Scientists, Splaiul Independen\c{t}ei 54, 050094
Bucharest, Romania

\textbf{E-Mail: cristinel.mortici@hotmail.com}%
\[
\]

\textbf{Abstract}
\end{center}

\begin{quotation}
{\small The aim of this work is to expose some asymptotic series associated
to some expressions involving the volume of the }$n${\small -dimensional
unit ball.}\newline
\textbf{All proofs and the methods used for improving the classical
inequalities announced in the final part of the first section are presented
in an extended form in a paper submitted by the author to a journal for
publication.}
\end{quotation}

\medskip



\section{Introduction and Motivation}

In the recent past, inequalities about the volume of the unit ball in $%
\mathbb{R}
^{n}:$%
\begin{equation}
\Omega _{n}=\frac{\pi ^{n/2}}{\Gamma \left( \frac{n}{2}+1\right) }\ \ \ \ \
\left( n\in 
\mathbb{N}
\right)   \label{ono}
\end{equation}%
have attracted the attention of many authors. See, \emph{e.g.}, \cite{alzer}-%
\cite{qiu}. Here $\Gamma $ denotes the Euler's gamma function defined for
every real number $x>0,$ by the formula:%
\[
\Gamma \left( x+1\right) =\int_{0}^{\infty }t^{x}e^{-t}dt,
\]%
while $%
\mathbb{N}
$ denotes the set of all positive integers and $%
\mathbb{N}
_{0}=%
\mathbb{N}
\cup \left\{ 0\right\} .$

Our products improve the following classical results:

\begin{itemize}
\item Chen and Lin \cite{chen} ($a=\frac{e}{2}-1,$ $b=\frac{1}{3})$:%
\[
\frac{1}{\sqrt{\pi \left( n+a\right) }}\left( \frac{2\pi e}{n}\right) ^{%
\frac{n}{2}}\leq \Omega _{n}<\frac{1}{\sqrt{\pi \left( n+b\right) }}\left( 
\frac{2\pi e}{n}\right) ^{\frac{n}{2}}\ \ \ \ \ \left( n\in 
\mathbb{N}
\right) ; 
\]

\item Borgwardt \cite{bor} ($a=0,$ $b=1),$ Alzer \cite{alzer2} and Qiu and
Vuorinen \cite{qiu} ($a=\frac{1}{2}$, $b=\frac{\pi }{2}-1$):%
\[
\sqrt{\frac{n+a}{2\pi }}\leq \frac{\Omega _{n-1}}{\Omega _{n}}\leq \sqrt{%
\frac{n+b}{2\pi }}\ \ \ \ \ \left( n\in 
\mathbb{N}
\right) ; 
\]

\item Alzer \cite{alzer2} ($\alpha ^{\ast }=\frac{3\pi \sqrt{2}}{4\pi +6},\
\beta ^{\ast }=\sqrt{2\pi }$):%
\[
\frac{\alpha ^{\ast }}{\sqrt{n}}\leq \frac{\Omega _{n}}{\Omega _{n-1}+\Omega
_{n+1}}<\frac{\beta ^{\ast }}{\sqrt{n}}\ \ \ \ \ \left( n\in 
\mathbb{N}
\right) ; 
\]

\item Chen and Lin \cite{chen} ($a=\frac{\pi \left( 1+\pi \right) ^{2}}{2}%
-1, $ $b=\frac{1}{2}+4\pi )$:%
\[
\sqrt{\frac{2\pi }{n+a}}\leq \frac{\Omega _{n}}{\Omega _{n-1}+\Omega _{n+1}}<%
\sqrt{\frac{2\pi }{n+b}}\ \ \ \ \ \left( n\in 
\mathbb{N}
\right) ; 
\]

\item Chen and Lin in \cite{chen} ($\lambda =1,$ $\mu =\frac{2\ln 2-\ln \pi 
}{2\ln 3-3\ln 2}$):%
\[
\left( 1+\frac{1}{2n}-\frac{3}{8n^{2}}\right) ^{\lambda }<\frac{\Omega
_{n}^{2}}{\Omega _{n-1}\Omega _{n+1}}\leq \left( 1+\frac{1}{2n}-\frac{3}{%
8n^{2}}\right) ^{\mu }\ \ \ \ \ \left( n\in 
\mathbb{N}
\right) ; 
\]

\item Anderson et al. \cite{anderson} and Klain and Rota \cite{klain}:%
\[
1<\frac{\Omega _{n}^{2}}{\Omega _{n-1}\Omega _{n+1}}<1+\frac{1}{n}\ \ \ \ \
\left( n\in 
\mathbb{N}
\right) ; 
\]

\item Alzer \cite{alzer2} ($\alpha =2-\log _{2}\pi ,$ $\beta =\frac{1}{2}$):%
\[
\left( 1+\frac{1}{n}\right) ^{\alpha }\leq \frac{\Omega _{n}^{2}}{\Omega
_{n-1}\Omega _{n+1}}<\left( 1+\frac{1}{n}\right) ^{\beta }\ \ \ \ \ \left(
n\in 
\mathbb{N}
\right) ; 
\]

\item Merkle \cite{merkle}:%
\[
\left( 1+\frac{1}{n+1}\right) ^{\frac{1}{2}}<\frac{\Omega _{n}^{2}}{\Omega
_{n-1}\Omega _{n+1}}\ \ \ \ \ \left( n\in 
\mathbb{N}
\right) ; 
\]

\item Chen and Lin \cite{chen} ($\alpha =\frac{1}{2},$ $\beta =\frac{2\ln
2-\ln \pi }{\ln 3-\ln 2}$):%
\[
\left( 1+\frac{1}{n+1}\right) ^{\alpha }<\frac{\Omega _{n}^{2}}{\Omega
_{n-1}\Omega _{n+1}}\leq \left( 1+\frac{1}{n+1}\right) ^{\beta }\ \ \ \ \
\left( n\in 
\mathbb{N}
\right) .
\]
\end{itemize}

\section{Classical results and new achievements}

\subsection{Asymptotic series and estimates for $\Omega _{n}$ and $\Omega
_{n}^{1/n}$}

Mortici \cite[Rel. 17]{m2} established the following asymptotic series as $%
n\rightarrow \infty :$%
\begin{eqnarray*}
\frac{1}{n}\ln \Omega _{n} &\sim &-\frac{n+1}{2n}\ln \frac{n}{2}+\frac{1}{2}%
\ln \left( \pi e\right) -\frac{\ln 2\pi }{2n} \\
&&-\left( \frac{1}{6n^{2}}-\frac{1}{45n^{4}}+\frac{8}{315n^{6}}-\frac{8}{%
105n^{8}}+\frac{1}{128n^{10}}+\cdots \right) .
\end{eqnarray*}%
The entire series is given below:

\begin{theorem}
The following asymptotic series holds true, as $n\rightarrow \infty :$%
\begin{equation}
\frac{1}{n}\ln \Omega _{n}\sim -\frac{n+1}{2n}\ln \frac{n}{2}+\frac{1}{2}\ln
\left( \pi e\right) -\frac{\ln 2\pi }{2n}-\sum_{j=1}^{\infty }\frac{%
2^{2j-2}B_{2j}}{j\left( 2j-1\right) n^{2j}}.  \label{py}
\end{equation}%
($B_{j}$ are the Bernoulli numbers).
\end{theorem}

The following double inequality \cite[Theorem 2]{m2} was presented$:$%
\begin{equation}
\alpha \left( n\right) <\frac{1}{n}\ln \Omega _{n}<\beta \left( n\right) \ \
\ \ \ \left( n\in \mathbb{%
\mathbb{N}
}\right) ,  \label{abba}
\end{equation}%
where%
\begin{eqnarray*}
\alpha \left( n\right)  &=&-\frac{n+1}{2n}\ln \frac{n}{2}+\frac{1}{2}\ln
\left( \pi e\right) -\frac{\ln 2\pi }{2n}-\lambda \left( n\right)  \\
\beta \left( n\right)  &=&-\frac{n+1}{2n}\ln \frac{n}{2}+\frac{1}{2}\ln
\left( \pi e\right) -\frac{\ln 2\pi }{2n}-\mu \left( n\right) ,
\end{eqnarray*}%
with%
\begin{eqnarray*}
\mu \left( n\right)  &=&\frac{1}{6n^{2}}-\frac{1}{45n^{4}}+\frac{8}{315n^{6}}%
-\frac{8}{105n^{8}} \\
\lambda \left( n\right)  &=&\mu \left( n\right) +\frac{1}{128n^{10}}.
\end{eqnarray*}

\begin{theorem}
The following inequality holds true$:$%
\[
\alpha \left( n\right) -\beta \left( n+1\right) >0\ \ \ \ \ \left( n\in 
\mathbb{%
\mathbb{N}
}\right) . 
\]%
In consequence, the sequence $\left\{ \Omega _{n}^{1/n}\right\} _{n\geq 1}$
decreases monotonically (to $0$).
\end{theorem}

\begin{theorem}
The following double inequality holds true, for every integer $n\geq 3$ in
the left-hand side and $n\geq 1$ in the right-hand side:%
\begin{equation}
\frac{1}{\sqrt{\pi \left( n+\theta \left( n\right) \right) }}\left( \frac{%
2\pi e}{n}\right) ^{\frac{n}{2}}<\Omega _{n}<\frac{1}{\sqrt{\pi \left( n+\nu
\left( n\right) \right) }}\left( \frac{2\pi e}{n}\right) ^{\frac{n}{2}},
\label{33}
\end{equation}%
where%
\begin{eqnarray*}
\theta \left( n\right)  &=&\frac{1}{3}+\frac{1}{18n}-\frac{31}{810n^{2}} \\
\nu \left( n\right)  &=&\theta \left( n\right) -\frac{139}{9720n^{3}}.
\end{eqnarray*}
\end{theorem}

Next we construct an asymptotic series for the ratio $\Omega
_{n}^{1/n}/\Omega _{n+1}^{1/\left( n+1\right) },$ then we give some lower
and upper bounds.

\begin{theorem}
The following asymptotic series holds true, as $n\rightarrow \infty :$%
\begin{equation}
\frac{1}{n}\ln \Omega _{n}-\frac{1}{n+1}\ln \Omega _{n+1}\sim \Psi \left(
n\right) -\sum_{j=1}^{\infty }\frac{\psi _{j}}{n^{j}}  \label{psi}
\end{equation}%
where%
\[
\Psi \left( n\right) =-\frac{n+1}{2n}\ln \frac{n}{2}+\frac{n+2}{2n+2}\ln 
\frac{n+1}{2}-\frac{\ln 2\pi }{2n\left( n+1\right) } 
\]%
and the coefficients $\psi _{j}$ are given by:%
\begin{eqnarray*}
\psi _{2t+1} &=&\sum_{k+2s=2t+1}\left( -1\right) ^{k+1}\left( 
\begin{array}{c}
k+2s-1 \\ 
k%
\end{array}%
\right) \ \ \ \ \ \left( t,k,s\in 
\mathbb{N}
\right) \\
\psi _{2t} &=&\frac{2^{2t-2}B_{2t}}{t\left( 2t-1\right) }+\sum_{k+2s=2t}%
\left( -1\right) ^{k+1}\left( 
\begin{array}{c}
k+2s-1 \\ 
k%
\end{array}%
\right) \ \ \ \ \ \left( t,k,s\in 
\mathbb{N}
\right) ,
\end{eqnarray*}%
with%
\[
\left( 
\begin{array}{c}
v \\ 
k%
\end{array}%
\right) =\frac{v\left( v-1\right) \cdots \left( v-k+1\right) }{k!}\ \ \ \ \
\left( v\in 
\mathbb{R}
,\text{ }k\in 
\mathbb{N}
_{0}\right) . 
\]
\end{theorem}

We have:%
\[
\frac{1}{n}\ln \Omega _{n}-\frac{1}{n+1}\ln \Omega _{n+1}\sim \Psi \left(
n\right) -\frac{1}{3n^{3}}+\frac{1}{2n^{4}}-\frac{26}{45n^{5}}+\frac{11}{%
18n^{6}}+\cdots .
\]

\begin{theorem}
The following double inequality holds true:%
\[
\Psi \left( n\right) -\frac{1}{3n^{3}}<\frac{1}{n}\ln \Omega _{n}-\frac{1}{%
n+1}\ln \Omega _{n+1}<\Psi \left( n\right) -\frac{1}{3n^{3}}+\frac{1}{2n^{4}}%
\ \ \ \ \ \left( n\in \mathbb{N}\right) . 
\]
\end{theorem}

\subsection{Asymptotic series and estimates for $\frac{\Omega _{n-1}}{\Omega
_{n}}$}

\begin{theorem}
The following asymptotic series holds true, as $n\rightarrow \infty :$%
\begin{equation}
\ln \frac{\Omega _{n-1}}{\Omega _{n}}=\frac{1}{2}\ln \frac{n}{2\pi }%
+\sum_{j=1}^{\infty }\frac{\mu _{j}}{n^{j}},  \label{miu}
\end{equation}%
where%
\[
\mu _{j}=\left( -1\right) ^{j}\left[ B_{j+1}\left( \frac{1}{2}\right)
-B_{j+1}\left( 1\right) \right] \frac{2^{j}}{j\left( j+1\right) }\ \ \ \ \
\left( j\in 
\mathbb{N}
\right) 
\]%
($B_{j}$ are the Bernoulli polynomials).
\end{theorem}

In a concrete form, (\ref{miu}) can be written as:%
\[
\ln \frac{\Omega _{n-1}}{\Omega _{n}}=\frac{1}{2}\ln \frac{n}{2\pi }+\left( 
\frac{1}{4n}-\frac{1}{24n^{3}}+\frac{1}{20n^{5}}-\frac{17}{112n^{7}}+\frac{31%
}{36n^{9}}+\cdots \right) . 
\]%
Remark that only odd powers of $n^{-1}$ appear in this series. This can be
justified using the following representation formulas of the Bernoulli
polynomials in terms of the Bernoulli numbers:%
\[
B_{t}\left( 1\right) =\left( -1\right) ^{t}B_{t}\ ,\ \ \ B_{t}\left( \frac{1%
}{2}\right) =\left( 2^{1-t}-1\right) B_{t}\ \ \ \ \ \left( t\in 
\mathbb{N}
\right) . 
\]%
See, \emph{e.g.}, \cite{wh}. As the Bernoulli numbers of odd order vanish,
it results that%
\[
B_{j+1}\left( \frac{1}{2}\right) =B_{j+1}\left( 1\right) =0\ \ \ \ \ \left(
j\in 2%
\mathbb{N}
\right) 
\]%
and consequently, $\mu _{j}=0,$ whenever $j$ is a positive even integer.

We present the following estimates:

\begin{theorem}
The following double inequality holds true:%
\begin{equation}
a\left( n\right) <\ln \frac{\Omega _{n-1}}{\Omega _{n}}<b\left( n\right) \ \
\ \ \ \left( n\in 
\mathbb{N}
\right) ,  \label{ab}
\end{equation}%
where%
\begin{eqnarray*}
a\left( n\right) &=&\frac{1}{2}\ln \frac{n}{2\pi }+\frac{1}{4n}-\frac{1}{%
24n^{3}}+\frac{1}{20n^{5}}-\frac{17}{112n^{7}} \\
b\left( n\right) &=&a\left( n\right) +\frac{31}{36n^{9}}.
\end{eqnarray*}
\end{theorem}

Further we deduce a new asymptotic series and one of the resulting double
inequality:

\begin{theorem}
The following asymptotic series holds true, as $n\rightarrow \infty :$%
\begin{equation}
\frac{\Omega _{n-1}}{\Omega _{n}}=\sqrt{\frac{n}{2\pi }}\left\{
\sum_{j=0}^{\infty }\frac{c_{j}}{n^{j}}\right\} ,  \label{cj}
\end{equation}%
where $c_{0}=1$ and%
\[
c_{j}=\frac{1}{j}\sum_{k=1}^{j}\left( -1\right) ^{k}\left[ B_{k+1}\left( 
\frac{1}{2}\right) -B_{k+1}\left( 1\right) \right] \frac{2^{k}}{k+1}c_{j-k}\
\ \ \ \ \left( j\in 
\mathbb{N}
\right) . 
\]
\end{theorem}

By listing the first terms, we get:%
\[
\frac{\Omega _{n-1}}{\Omega _{n}}=\sqrt{\frac{n}{2\pi }}\left( 1+\frac{1}{4n}%
+\frac{1}{32n^{2}}-\frac{5}{128n^{3}}-\frac{21}{2048n^{4}}+\frac{399}{%
8192n^{5}}+\frac{869}{65\,536n^{6}}+\cdots \right) . 
\]%
Related to this asymptotic expansion, we prove the following estimates:

\begin{theorem}
The following double inequality holds true, for every integer $n\geq 12$ in
the left-hand side and $n\geq 1$ in the right-hand side:%
\begin{equation}
\sqrt{\frac{n}{2\pi }}c\left( n\right) <\frac{\Omega _{n-1}}{\Omega _{n}}<%
\sqrt{\frac{n}{2\pi }}d\left( n\right) ,  \label{cd}
\end{equation}%
where%
\begin{eqnarray*}
c\left( n\right) &=&1+\frac{1}{4n}+\frac{1}{32n^{2}}-\frac{5}{128n^{3}}-%
\frac{21}{2048n^{4}}+\frac{399}{8192n^{5}} \\
d\left( n\right) &=&c\left( n\right) +\frac{869}{65\,536n^{6}}.
\end{eqnarray*}
\end{theorem}

\begin{theorem}
The following asymptotic formula holds true as $n\rightarrow \infty :$%
\begin{equation}
\frac{\Omega _{n-1}}{\Omega _{n}}=\sqrt{\frac{n+\frac{1}{2}}{2\pi }%
+\sum_{j=1}^{\infty }\frac{s_{j}}{n^{j}}},
\end{equation}%
where%
\[
s_{j}=\frac{1}{2\pi }\sum_{k=0}^{j+1}c_{k}c_{j+1-k}\ \ \ \ \ \left( j\in 
\mathbb{N}
\right) .
\]
\end{theorem}

The first terms are indicated below:%
\begin{eqnarray*}
\frac{\Omega _{n-1}}{\Omega _{n}} &=&\left( \frac{n+\frac{1}{2}}{2\pi }+%
\frac{1}{16\pi n}-\frac{1}{32\pi n^{2}}-\frac{5}{256\pi n^{3}}\right. \\
&&\left. +\frac{23}{512\pi n^{4}}+\frac{53}{2048\pi n^{5}}-\allowbreak \frac{%
593}{4096\pi n^{6}}-\frac{5165}{65\,536\pi n^{7}}+\cdots \right) ^{\frac{1}{2%
}}.
\end{eqnarray*}%
By truncation of this series, increasingly accurate under- and upper-
approximations for the ratio $\frac{\Omega _{n-1}}{\Omega _{n}}$ are
obtained. As an example, we show the following:

\begin{theorem}
The following double inequality holds true, for every integer $n\geq 1$ in
the left-hand side and $n\geq 2$ in the right-hand side:%
\[
\sqrt{\frac{n+\frac{1}{2}}{2\pi }+\frac{1}{16\pi n}-\frac{1}{32\pi n^{2}}-%
\frac{5}{256\pi n^{3}}}<\frac{\Omega _{n-1}}{\Omega _{n}}<\sqrt{\frac{n+%
\frac{1}{2}}{2\pi }+\frac{1}{16\pi n}-\frac{1}{32\pi n^{2}}}. 
\]
\end{theorem}

\begin{theorem}
The following double inequality holds true$:$%
\begin{equation}
\sqrt{\frac{2\pi }{n+4\pi +\frac{1}{2}}+\varepsilon _{1}\left( n\right) }<%
\frac{\Omega _{n}}{\Omega _{n-1}+\Omega _{n+1}}<\sqrt{\frac{2\pi }{n+4\pi +%
\frac{1}{2}}+\varepsilon _{2}\left( n\right) }\ \ \ \ \ \left( n\in 
\mathbb{N}
\right) ,  \label{o6}
\end{equation}%
where%
\begin{eqnarray*}
\varepsilon _{1}\left( n\right) &=&-\frac{\frac{1}{4}\pi -4\pi ^{2}+8\pi ^{3}%
}{n^{3}} \\
\varepsilon _{2}\left( n\right) &=&\varepsilon _{1}\left( n\right) +\frac{%
\frac{3}{8}\pi -7\pi ^{2}-12\pi ^{3}+64\pi ^{4}}{n^{4}}.
\end{eqnarray*}
\end{theorem}

\subsection{Asymptotic series and estimates for $\frac{\Omega _{n}^{2}}{%
\Omega _{n-1}\Omega _{n+1}}$}

We start this section by establishing new asymptotic expansions for the
ratio $\frac{\Omega _{n}^{2}}{\Omega _{n-1}\Omega _{n+1}}$\ and some
associated inequalities.

\begin{theorem}
The following asymptotic series holds true, as $n\rightarrow \infty :$%
\begin{equation}
\ln \frac{\Omega _{n}^{2}}{\Omega _{n-1}\Omega _{n+1}}=\sum_{j=1}^{\infty }%
\frac{\lambda _{j}}{n^{j}},  \label{laj}
\end{equation}%
where%
\[
\lambda _{j}=\left( -1\right) ^{j}\left\{ 2B_{j+1}\left( 1\right)
-B_{j+1}\left( \frac{1}{2}\right) -B_{j+1}\left( \frac{3}{2}\right) \right\} 
\frac{2^{j}}{j\left( j+1\right) }\ \ \ \ \ \left( j\in 
\mathbb{N}
\right) . 
\]
\end{theorem}

As the first terms in this series are%
\[
\ln \frac{\Omega _{n}^{2}}{\Omega _{n-1}\Omega _{n+1}}=\frac{1}{2n}-\frac{1}{%
2n^{2}}+\frac{5}{12n^{3}}-\frac{1}{4n^{4}}+\frac{1}{10n^{5}}-\frac{1}{6n^{6}}%
+\cdots , 
\]%
we are entitled to present the following estimates:

\begin{theorem}
The following double inequality holds true$:$%
\begin{equation}
p\left( n\right) <\ln \frac{\Omega _{n}^{2}}{\Omega _{n-1}\Omega _{n+1}}%
<q\left( n\right) \ \ \ \ \ \left( n\in 
\mathbb{N}
\right) ,  \label{pq}
\end{equation}%
where%
\begin{eqnarray*}
p\left( n\right) &=&\frac{1}{2n}-\frac{1}{2n^{2}}+\frac{5}{12n^{3}}-\frac{1}{%
4n^{4}}+\frac{1}{10n^{5}}-\frac{1}{6n^{6}} \\
q\left( n\right) &=&p\left( n\right) +\frac{1}{6n^{6}}.
\end{eqnarray*}
\end{theorem}

\begin{theorem}
The following asymptotic series holds true, as $n\rightarrow \infty :$%
\[
\frac{\Omega _{n}^{2}}{\Omega _{n-1}\Omega _{n+1}}=\sum_{j=0}^{\infty }\frac{%
d_{j}}{n^{j}}, 
\]%
where $d_{0}=1$ and $d_{j}^{\prime }$s, $j\in 
\mathbb{N}
,$ are defined by the recursive relation:%
\[
d_{j}=\frac{1}{j}\sum_{k=1}^{j}\left( -1\right) ^{k}\left[ 2B_{k+1}\left(
1\right) -B_{k+1}\left( \frac{1}{2}\right) -B_{k+1}\left( \frac{3}{2}\right) %
\right] \frac{2^{k}}{k+1}d_{j-k}. 
\]
\end{theorem}

More exactly, we have:%
\begin{equation}
\frac{\Omega _{n}^{2}}{\Omega _{n-1}\Omega _{n+1}}=1+\frac{1}{2n}-\frac{3}{%
8n^{2}}+\frac{3}{16n^{3}}+\frac{3}{128n^{4}}+\cdots   \label{pop}
\end{equation}

We propose the following estimates associated to the series (\ref{pop}):

\begin{theorem}
The following double inequality holds true, for every integer $n\geq 6$ in
the left-hand side and $n\geq 1$ in the right-hand side:%
\begin{equation}
r\left( n\right) <\frac{\Omega _{n}^{2}}{\Omega _{n-1}\Omega _{n+1}}<s\left(
n\right) ,  \label{rs}
\end{equation}%
where%
\begin{eqnarray*}
r\left( n\right) &=&1+\frac{1}{2n}-\frac{3}{8n^{2}}+\frac{3}{16n^{3}} \\
s\left( n\right) &=&r\left( n\right) +\frac{3}{128n^{4}}.
\end{eqnarray*}
\end{theorem}

\begin{theorem}
The following asymptotic formula holds true as $n\rightarrow \infty :$%
\begin{equation}
\frac{\Omega _{n}^{2}}{\Omega _{n-1}\Omega _{n+1}}\sim \left( 1+\frac{1}{n}%
\right) ^{\sum_{j=0}^{\infty }\frac{t_{j}}{n^{j}}},  \label{teen}
\end{equation}%
where $t_{0}=\frac{1}{2}$ and $t_{j}^{\prime }$s, $j\in 
\mathbb{N}
,$ are the solution of the infinite system:%
\[
\sum_{j=1}^{m}\left( -1\right) ^{j+1}\frac{t_{m-j}}{j}=\lambda _{m}\ \ \ \ \
\left( m\in 
\mathbb{N}
\right) .
\]
\end{theorem}

\begin{theorem}
The following double inequality holds true, for every integer $n\geq 5$ in
the left-hand side and $n\geq 1$ in the right-hand side:%
\begin{equation}
\left( 1+\frac{1}{n}\right) ^{\frac{1}{2}-\frac{1}{4n}+\frac{1}{8n^{2}}}<%
\frac{\Omega _{n}^{2}}{\Omega _{n-1}\Omega _{n+1}}<\left( 1+\frac{1}{n}%
\right) ^{\frac{1}{2}-\frac{1}{4n}+\frac{1}{8n^{2}}+\frac{1}{48n^{3}}}.
\label{oo}
\end{equation}
\end{theorem}

Here we only list the following results:%
\[
\frac{\Omega _{n}^{2}}{\Omega _{n-1}\Omega _{n+1}}\sim \left( 1+\frac{1}{n+1}%
\right) ^{\frac{1}{2}+\frac{1}{4n}-\frac{3}{8n^{2}}+\frac{23}{48n^{3}}-\frac{%
15}{32n^{4}}+\cdots }
\]%
and%
\[
\frac{\Omega _{n}^{2}}{\Omega _{n-1}\Omega _{n+1}}\sim \left( 1+\frac{1}{n+1}%
\right) ^{\frac{1}{2}+\frac{1}{4\left( n+1\right) }-\frac{1}{8\left(
n+1\right) ^{2}}-\frac{1}{48\left( n+1\right) ^{3}}+\frac{3}{32\left(
n+1\right) ^{4}}+\cdots }.
\]

\textbf{Acknowledgements. 1. }This work was supported by a Grant of the
Romanian National Authority for Scientific Research CNCS-UEFISCDI, no.
PN-II-ID-PCE-2011-3-0087. The computations made in this paper were performed
using the Maple software for symbolic computation.

2. \textbf{All proofs and the methods used for improving the classical
inequalities announced in the final part of the first section are presented
in an extended form in a paper submitted by the author to a journal for
publication.}

\end{document}